\documentclass[11pt,reqno,twoside]{amsart}

\usepackage{amsmath, amssymb, amscd, amsthm}

\theoremstyle{plain}
\newtheorem{thm}{Theorem}[section]
\newtheorem{prop}[thm]{Proposition}

\theoremstyle{remark}
\newtheorem{rem}[thm]{Remark}
\newcounter{sspar}[subsection]
\renewcommand\thesspar{(\thesubsection.\arabic{sspar})}
    {\par\ \newline
     \vskip-\baselineskip\vskip.1truecm
     \noindent\refstepcounter{sspar}
     \noindent\textbf{\thesspar} \ignorespaces}
    {\vskip-\baselineskip
    \ignorespaces}
    {\refstepcounter{sspar}
     \textup{\textbf{\thesspar}} \ignorespaces}
    {\vskip-\baselineskip
    \ignorespaces}
\newcommand{\R}{{\mathbb R}}
\newcommand{\N}{{\mathbb N}}
\newcommand{\C}{{\mathbb C}}
\newcommand{\Z}{{\mathbb Z}}
\newcommand{\al}{\alpha}
\newcommand{\be}{\beta}

\newcommand{\Ga}{\Gamma}

\newcommand{\De}{\Delta}

\newcommand{\eps}{\varepsilon}
\newcommand{\si}{\sigma}

\newcommand{\la}{\lambda}
\newcommand{\La}{\Lambda}

\newcommand{\Om}{\Omega}
\newcommand{\lt}{\ell^2}
\newcommand{\F}[5]{\,_{#1}F_{#2} \left( \genfrac{.}{.}{0pt}{}{#3}{#4}
\ ;#5 \right)}
\newcommand{\inprod}[2]{\langle #1,#2 \rangle}
\newcommand{\oz}{\overline{z}}
\newcommand{\M}{\hat{M}}
\newcommand{\m}{\hat{m}}

\newcommand{\tensor}{\otimes}
\newcommand{\Sh}{\hat{S}}
\newcommand{\ow}{\overline{w}}
\newcommand{\w}{\widetilde{w}}
\newcommand{\dirint}{\sideset{}{^\oplus}\int\limits_{0\ }^{\ \infty}}

\newcommand{\su}{\mathfrak{su}}
\DeclareMathOperator{\supp}{supp}

\numberwithin{equation}{section}

\title{Meixner functions and polynomials related to Lie algebra
representations}
\author{Wolter Groenevelt}

\author{Erik Koelink}

\address{Technische Universiteit Delft, ITS-TWA \\
Postbus 5031, 2600 GA Delft, The Netherlands}
\email{W.G.M.Groenevelt@its.tudelft.nl  \\ H.T.Koelink@its.tudelft.nl}

\date{September 25, 2001}

\begin{document}
\begin{abstract}
The decomposition of the tensor product of a positive and a negative
discrete series representation of the Lie algebra $\su(1,1)$ is a
direct integral over the principal unitary series representations. In
the decomposition discrete terms can occur, and the discrete terms are
a finite number of discrete series representations or one complementary
series representation. The interpretation of Meixner functions and
polynomials as overlap coefficients in the four classes of
representations and the Clebsch-Gordan decomposition, lead to a general
bilinear generating function for the Meixner polynomials. Finally,
realizing the positive and negative discrete series representations as
operators on the spaces of holomorphic and anti-holomorphic functions
respectively, a non-symmetric type Poisson kernel is found for the
Meixner functions.
\end{abstract}

\maketitle

\section{Introduction}
The representation theory of Lie algebras is intimately related to
special functions of hypergeometric type, see e.g. Vilenkin and Klimyk
\cite{VK} and references therein. In this paper we consider the Lie
algebra $\su(1,1)$ spanned by $H$, $B$ and $C$. In \cite[\S3]{KJ1} Van
der Jeugt and the second author considered the self-adjoint element
$X_a = -2a H+B-C$, $a \geq 0$, in the double and triple tensor product
representation of positive discrete series representations, to find
convolution identities for orthogonal polynomials, see also \cite{Jeu}
where the case $a=1$ is treated. The method used is based on an idea of
Granovskii and Zhedanov \cite{GZ}; (generalized) eigenvectors of a
certain element of the Lie algebra, which acts as a three-term
recurrence operator in an irreducible representation of the Lie
algebra, can be used to obtain identities for special functions. In
this paper we apply this idea to the tensor product of a positive and a
negative discrete series representation of $\su(1,1)$.

There are four classes of irreducible unitary representations of
$\su(1,1)$; positive and negative discrete series, principal unitary
series and complementary series representations. In the discrete series
representations, the overlap coefficients between the standard basis
and the basis of eigenvectors of $X_c$ are all polynomials; Meixner,
Meixner-Pollaczek and Laguerre polynomials, corresponding to $a>1$,
$a<1$ and $a=1$ respectively (see Masson and Repka \cite{MR}). This is
no longer the case for the principal unitary and the complementary
series. Since both these representations appear in the decomposition of
the tensor product of a positive and a negative discrete series
representation, we have to use non-polynomial functions. In this paper
we concentrate on the case $a>1$. Then we need Meixner functions, which
are found in  \cite{Koe}, \cite{MR}, using spectral theory of Jacobi
operators.

In section 2 we decompose the tensor product of a positive and a
negative discrete series representation into a direct integral over the
principal series representations. In the decomposition discrete terms
can occur, and the discrete terms are a finite number of discrete
series representations or one complementary series representation. To
find the decomposition, we consider the spectrum of the Casimir
operator in the tensor product. As a consequence we find that the
continuous dual Hahn polynomials have a natural interpretation as
Clebsch-Gordan coefficients for $\su(1,1)$.

In section 3 we determine the eigenvectors of $X_a$ in the four
different representations, as well as the eigenvector of $X_a$ in the
tensor product representation. Then we determine the Clebsch-Gordan
coefficients for these eigenvectors, which turn out to be continuous
dual Hahn polynomials. Finally we obtain a bilinear generating function
for Meixner polynomials.

In section 4 use the holomorphic and anti-holomorphic realization of
the positive and negative discrete series representations. This
realization was also used in \cite{JJ} to find generating functions
from a tensor product of two positive discrete series representations.
Now the eigenvectors of $X_a$ become known generating functions for
Meixner polynomials. Using the Clebsch-Gordan decomposition from
section 2, we obtain a non-symmetric type Poisson kernel for the
Meixner functions.

\textit{Acknowledgment.} We thank Hjalmar Rosengren for useful
discussions.

\section{Decomposition of tensor product representations}
In this section we consider the Lie algebra $\su(1,1)$ and the four
classes of irreducible unitary representations. Then we decompose the
tensor product of a positive and a negative discrete series
representation into a direct integral of principal unitary series. In
the decomposition discrete terms can occur, and the discrete terms are
a finite number of discrete series representations or one complementary
series representation. We also determine the Clebsch-Gordan
coefficients, which turn out to be continuous dual Hahn polynomials.

\subsection{The Lie algebra $\su(1,1)$}
The Lie algebra $\su(1,1)$ is given by
\[
[H,B]=2B, \quad [H,C]=-2C, \quad [B,C]=H.
\]
The $\ast$-structure is defined by $H^*=H$ and $B^*=-C$. The Casimir
operator $\Om$ is a central self-adjoint element of the universal
enveloping algebra $U\big(\su(1,1)\big)$ and
\begin{equation} \label{Casimir}
 \Omega = -\frac{1}{4}(H^2+2H+4CB) = -\frac{1}{4}(H^2-2H+4BC).
\end{equation}

There are four classes of irreducible unitary representations of
$\su(1,1)$, see \cite[\S6.4]{VK}. The positive discrete series
representations $\pi_k^+$ are representations labeled by $k>0$. The
representation space is $\lt(\Z_{\geq 0})$ with orthonormal basis
$\{e_n\}_{n\in \Z_{\geq 0}}$. The action is given by
\begin{equation} \label{pos}
\begin{split}
\pi_k^+(H) e_n =&\ 2(k+n) e_n,  \\ \pi_k^+(B) e_n =&\
\sqrt{(n+1)(2k+n)} e_{n+1},  \\ \pi_k^+(C) e_n =&\ -\sqrt{n(2k+n-1)}
e_{n-1},  \\ \pi_k^+(\Omega) e_n =&\ k(1-k) e_n.
\end{split}
\end{equation}

The negative discrete series representations $\pi_{k}^-$ are labeled by
$k>0$. The representation space is $\lt(\Z_{\geq 0})$ with orthonormal
basis $\{e_n\}_{n \in \Z_{\geq 0}}$. The action is given by
\begin{equation} \label{neg}
\begin{split}
\pi_{k}^-(H) e_n =&\ -2(k+n) e_n, \\ \pi_{k}^-(B) e_n =&\
-\sqrt{n(2k+n-1)} e_{n-1},  \\ \pi_{k}^-(C) e_n =&\ \sqrt{(n+1)(2k+n)}
e_{n+1},  \\ \pi_{k}^-(\Omega) e_n =&\ k(1-k) e_n.
\end{split}
\end{equation}

The principal series representations $\pi^{\rho,\eps}$ are labeled by
$\eps \in [0,1)$ and $\rho \geq 0$, where $(\rho,\eps) \neq
(0,\frac{1}{2})$. The representation space is $\lt(\Z)$ with
orthonormal basis $\{e_n\}_{n \in \Z}$. The action is given by
\begin{equation} \label{princ}
\begin{split}
\pi^{\rho, \eps}(H) e_n =&\ 2(\eps+n) e_n, \\ \pi^{\rho, \eps}(B) e_n
=&\ \sqrt{(n+\eps+\frac{1}{2}-i\rho) (n+\eps+\frac{1}{2}+i\rho)}
e_{n+1},
\\ \pi^{\rho, \eps}(C) e_n =&\
-\sqrt{(n+\eps-\frac{1}{2}-i\rho) (n+\eps-\frac{1}{2}+i\rho)} e_{n-1},
\\ \pi^{\rho, \eps}(\Omega) e_n =&\ (\rho^2+\frac{1}{4}) e_n.
\end{split}
\end{equation}

The complementary series representations $\pi^{\lambda,\eps}$ are
labeled by $\eps$ and $\lambda$, where $\eps \in [0,\frac{1}{2})$ and
$\lambda \in (-\frac{1}{2},-\eps)$  or $\eps \in (\frac{1}{2},1)$ and
$\lambda \in (-\frac{1}{2},\eps-1)$. The representation space is
$\lt(\Z)$ with orthonormal basis $\{e_n\}_{n \in \Z}$. The action is
given by
\begin{equation} \label{comp}
\begin{split}
\pi^{\lambda,\eps}(H) e_n =&\ 2(\eps+n) e_n, \\ \pi^{\lambda,\eps}(B)
e_n =&\ \sqrt{(n+\eps+1+\lambda) (n+\eps-\lambda)} e_{n+1}, \\
\pi^{\lambda,\eps}(C) e_n =&\ -\sqrt{(n+\eps+\lambda)
(n+\eps-\lambda-1)} e_{n-1},  \\ \pi^{\lambda,\eps}(\Omega) e_n =&\
-\lambda(1+\lambda) e_n.
\end{split}
\end{equation}
Notice that for $\la = -\frac{1}{2}+i\rho$ the actions in the principal
series and in the complementary series are the same.

The operators are unbounded, defined on the space of finite linear
combinations of the basisvectors. From general theory, see Schm\" udgen
\cite[Ch.8, Ch.10]{Sch}, it follows that all operators are closable.

\subsection{Tensor products of positive and negative discrete series
representations} In order to decompose the tensor product of a positive
and a negative discrete series representation we need the continuous
dual Hahn polynomials, which are a special case of the Wilson
polynomials, see Wilson \cite[\S 4]{Wil}; see also \cite[\S 6.10]{AAR},
\cite{KS};
\begin{equation} \label{cont dHahn}
S_n(y;a,b,c) = (a+b)_n (a+c)_n\F{3}{2}{-n,a+ix,a-ix}{a+b,a+c}{1},
\quad x^2=y.
\end{equation}
For real parameters $a$, $b$, $c$, with $a+b$, $a+c$, $b+c$ positive,
the continuous dual Hahn polynomials are orthogonal with respect to a
positive measure, supported on a subset of $\R$.\\ The orthonormal
continuous dual Hahn polynomials
\[
\Sh_n(y) = \Sh_n(y;a,b,c) = \frac{(-1)^n S_n(y;a,b,c)} {\sqrt{n!(a+b)_n
(a+c)_n (b+c)_n}}
\]
satisfy the three-term recurrence relation
\begin{equation} \label{S rec}
y \Sh_n(y) = a_n \Sh_{n+1}(y) +b_n \Sh_n(y) +a_{n-1}\Sh_{n-1}(y),
\end{equation}
where
\begin{eqnarray*}
a_n & = & \sqrt{(n+1)(n+a+b)(n+a+c)(n+b+c)}, \\ b_n & = &
2n^2+2n(a+b+c-\frac{1}{2})+ab+ac+bc.
\end{eqnarray*}
By Kummer's transformation, see e.g. \cite[Cor. 3.3.5]{AAR}, the
polynomials $S_n$ and $\Sh_n$ are symmetric in $a$, $b$ and $c$.

Without loss of generality we assume that $a$ is the smallest of the
real parameters $a$, $b$ and $c$. Let $d\mu^2(\cdot;a,b,c)$ be the
measure defined by
\begin{align*}
\int_{\R} f(y) d&\mu^2(y;a,b,c)  =
 \frac{1}{2\pi}\int_0^\infty \left|
\frac{\Gamma(a+ix)\Gamma(b+ix)\Gamma(c+ix)}{\Gamma(2ix)} \right|^2
\frac{f(x^2)}{\Gamma(a+b) \Gamma(a+c) \Gamma(b+c)}\,d x\\ + &
\frac{\Gamma(b-a) \Gamma(c-a)}{\Gamma(-2a) \Gamma(b+c)} \sum_{k=0}^K
(-1)^k \frac{(2a)_k (a+1)_k (a+b)_k (a+c)_k} {(a)_k (a-b+1)_k (a-c+1)_k
k!} f(-(a+k)^2),
\end{align*}
where $K$ is the largest non-negative integer such that $a+K<0$. In
particular, the measure $d\mu^2(\cdot;a,b,c)$ is absolutely continuous
if $a\geq0$. The measure is positive under the conditions $a+b>0$,
$a+c>0$ and $b+c>0$. Then the polynomials $\Sh_n(y;a,b,c)$ are
orthonormal with respect to the measure $d\mu^2(y;a,b,c)$.\\ Later on
we also use the notation $d \mu(\cdot;a,b,c)$, by this we denote the
measure defined by
\begin{align*}
\int_{\R} f(y) d&\mu(y;a,b,c)  =
 \frac{1}{\sqrt{2\pi}}\int_0^\infty \left|
\frac{\Gamma(a+ix)\Gamma(b+ix)\Gamma(c+ix)} {\Gamma(2ix)} \right|
\frac{f(x^2)} {\sqrt{\Gamma(a+b) \Gamma(a+c) \Gamma(b+c)}} \,d x \\ &+
\sqrt{ \frac{\Gamma(b-a) \Gamma(c-a)}{\Gamma(-2a) \Gamma(b+c)} }
\sum_{k=0}^K \sqrt{ (-1)^k \frac{(2a)_k (a+1)_k (a+b)_k (a+c)_k} {(a)_k
(a-b+1)_k (a-c+1)_k k!} }\  f(-(a+k)^2).
\end{align*}
Since $d\mu^2(\cdot;a,b,c)$ is a positive measure for $a+b>0$, $a+c>0$
and $b+c>0$, the measure $d \mu (\cdot;a,b,c)$ is a positive measure.\\

In order to determine the decomposition of the tensor product
representation $\pi_{k_1}^+ \tensor \pi_{k_2}^-$, we calculate the
action of the Casimir operator $\Om$. In the tensor product, we need
the comultiplication $\De$. Recall that for an element $Y$ of the Lie
algebra $\su(1,1)$, we have $\Delta(Y)=1 \tensor Y + Y \tensor 1$ and
$\De$ is extended to $U\big(\su(1,1)\big)$ as an algebra homomorphism.
Then from \eqref{Casimir} we obtain
\[
\Delta(\Omega) = 1 \tensor \Omega + \Omega \tensor 1 - \frac{1}{2} H
\tensor H -( C \tensor B + B \tensor C).
\]
We use this expression, \eqref{pos} and $\eqref{neg}$, to calculate the
action of $\Om$ in the tensor product
\begin{align*}
\pi_{k_1}^+ \tensor\pi_{k_2}^- \big(\Delta(\Omega) \big)&\, e_n \tensor
e_m
=
\\ -&\sqrt{(n+1)(2k_1+n)(m+1)(2k_2+m)}\, e_{n+1} \tensor e_{m+1}  \\ +&
\big(k_1(1-k_1) + k_2(1-k_2) + 2(k_1+n)(k_2+m) \big)\, e_n \tensor e_m
\\-& \sqrt{n(2k_1+n-1)\,m(2k_2+m-1)}\, e_{n-1} \tensor e_{m-1},
\end{align*}
We see that $\pi_{k_1}^+ \tensor\pi_{k_2}^- \big(\Delta(\Omega) \big)$
is a three-term recurrence operator, which we can identify with the
three-term recurrence relation for the continuous dual Hahn
polynomials. In order to do so we define for $p \in \Z$ and $n \in
\Z_{\geq 0}$, elements $f_n^p
\in
 \lt(\Z_{\geq 0}) \tensor \lt(\Z_{\geq 0})$ by
\[
f_n^p =
\begin{cases}
(-1)^n e_{n} \tensor e_{n-p}\,, & p \leq 0, \\ (-1)^n e_{n+p} \tensor
e_{n}\,, & p \geq0.
\end{cases}
\]
Note that $\{f_n^p\}_{n \in \Z_{\geq 0},p \in \Z}$ is an orthonormal
basis for $\lt(\Z_{\geq 0}) \tensor \lt(\Z_{\geq 0})$. For fixed $p$
define the space
\[
H_p = \overline{ \C \{f_n^p \ | \ n \in \Z_{\geq 0}\} },
\]
then $\pi_{k_1}^+ \tensor\pi_{k_2}^- \big(\Delta(\Omega) \big)$ is an
unbounded symmetric tridiagonal operator on $H_p$ with domain the space
of finite linear combinations of $f_n^p$.

We define operators $\La_p: H_p \rightarrow L^2(\R, d\mu^2(y;p) )$ by
\[
f_n^p \mapsto
\begin{cases}
\Sh_n(y;k_1-k_2+\frac{1}{2},k_1+k_2-\frac{1}{2},
k_2-k_1-p+\frac{1}{2}), & p \leq 0, \\
\Sh_n(y;k_2-k_1+\frac{1}{2},k_1+k_2-\frac{1}{2},
k_1-k_2+p+\frac{1}{2}), & p \geq 0,
\end{cases}
\]
where $d\mu^2(y;p)$ denotes the measure $d\mu^2(y;a,b,c)$ with
\begin{equation} \label{param}
a =
\begin{cases}
k_1-k_2+\frac{1}{2}, & p \leq 0, \\ k_2-k_1+\frac{1}{2}, & p \geq 0,
\end{cases} \quad b = k_1+k_2-\frac{1}{2}, \quad
c =
\begin{cases}
k_2-k_1-p+\frac{1}{2}, & p \leq 0, \\ k_1-k_2+p+\frac{1}{2}, & p \geq
0.
\end{cases}
\end{equation}
\begin{prop} \label{La p}
For every $p \in \Z$ the operator $\La_p$ is unitary and intertwines
$\pi_{k_1}^+ \tensor\pi_{k_2}^- \big(\Delta(\Omega) \big)$ acting on
$H_p$, with $M_{y+\frac{1}{4}}$ acting on $L^2(\R, d\mu^2(y;p))$.
\end{prop}
Here $M_g$ denotes multiplication by the function $g$, so $M_g f(y)=
g(y)f(y)$.
\begin{proof}
Letting $\pi_{k_1}^+ \tensor\pi_{k_2}^- \big(\Delta(\Omega) \big)$ act
on an element $f_n^p$, we see that $\pi_{k_1}^+ \tensor\pi_{k_2}^-
\big(\Delta(\Omega) \big)\big|_{H_p}$ is an unbounded Jacobi operator,
see Akhiezer \cite{Akh}. Comparing the coefficients with the three-term
recurrence relation \eqref{S rec} for the continuous dual Hahn
polynomials, we see that we can identify $\pi_{k_1}^+
\tensor\pi_{k_2}^- \big(\Delta(\Omega) \big)|_{H_p}-\frac{1}{4}$ with
\eqref{S rec} using \eqref{param} by a straightforward computation.
Since the corresponding moment problem is determinate, $\pi_{k_1}^+
\tensor\pi_{k_2}^- \big(\Delta(\Omega) \big)\big|_{H_p}$ is essentially
self-adjoint and the operator $\La_p$ intertwines $\pi_{k_1}^+
\tensor\pi_{k_2}^- \big(\Delta(\Omega) \big)\big|_{H_p}$ with
$M_{y+\frac{1}{4}}$. Since the continuous dual Hahn polynomials $\Sh_n$
are dense in  $L^2(\R, d\mu^2(y;p) )$, $\La_p$ maps one orthonormal
basis onto another. Therefore $\La_p$ is a unitary operator.
\end{proof}

Using proposition \ref{La p} we can find the spectrum of $\pi_{k_1}^+
\tensor\pi_{k_2}^- \big(\Delta(\Omega) \big)$, which is an essentially
self-adjoint operator on $\lt(\Z_{\geq 0}) \tensor \lt(\Z_{\geq 0})$.
From the orthonormality measure for the continuous dual Hahn
polynomials, we see that the spectrum has a continuous part
$[\frac{1}{4}, \infty)$ and a finite discrete part, depending on the
parameters $k_1$ and $k_2$. From now on we assume that $k_1 \leq k_2$.

For the continuous part of the spectrum, we recognize the action of
$\Om$ in the principal series representation, using $y = \rho^2$. From
the action of $H$ in the tensor product, we obtain $\eps = k_1-k_2+L$,
where $L$ is the unique non-negative integer such that $\eps \in
[0,1)$.

 If $k_1+k_2 < \frac{1}{2}$ the point spectrum contains one
term: $\frac{1}{4}-(k_1+k_2-\frac{1}{2})^2$. Using $\la = -k_1-k_2$ we
recognize here the action of $\Om$ in the complementary series
representation labeled by $\la$ and $\eps$.

For $k_1-k_2<-\frac{1}{2}$ the spectrum has $K+1$ discrete terms for $p
\leq 0$ and $K+1-p$ terms for $0 \leq p \leq K$. If $p >K$ there is no
discrete part. Here $K$ is the largest integer such that
$k_1-k_2+\frac{1}{2}+K<0$. All the terms in the discrete part of the
spectrum have the form $\frac{1}{4}-(k_1-k_2+\frac{1}{2}+j)^2$,
$j=0,1,\ldots,K$. From this we recognize the action of $\Omega$ in a
discrete series representation $\pi_{k_2-k_1-j}$. From the action of
$H$, we find that this is a negative discrete series representation.

Now we can prove the following theorem.

\begin{thm} \label{decomp}
For $k_1 \leq k_2$ the decomposition of the tensor product of positive
and negative discrete series representations of $\su(1,1)$ is
\begin{align*}
\pi_{k_1}^+ \tensor\pi_{k_2}^- &\cong \dirint \pi^{\rho,\eps} d \rho, &
k_1-k_2 \geq -\frac{1}{2}, k_1+k_2 \geq \frac{1}{2}, \\ \pi_{k_1}^+
\tensor\pi_{k_2}^- &\cong \dirint \pi^{\rho,\eps} d \rho \oplus
\pi^{\la, \eps}, & k_1+k_2<\frac{1}{2}, \\ \pi_{k_1}^+
\tensor\pi_{k_2}^- &\cong \dirint \pi^{\rho,\eps} d \rho \oplus
\bigoplus_{\substack{j \in \Z_{\geq 0}\\ k_2-k_1-\frac{1}{2}-j>0}}
\pi_{k_2-k_1-j}^-, & k_1-k_2
< -\frac{1}{2},
\end{align*}
where $\eps = k_1-k_2+L$, $L$ is the unique integer such that $\eps \in
[0,1)$ and $\la = -k_1-k_2$. Further, under the identification above,
\begin{equation} \label{decomp vec}
e_{n_1} \tensor e_{n_2} = (-1)^{n_2} \int_\R \Sh_n(y;n_1-n_2)
e_{n_1-n_2-L}\, d\mu(y;n_1-n_2),
\end{equation}
where $n= \min\{n_1,n_2\}$, $\Sh_n(y;p)$ is a continuous dual Hahn
polynomial,
\[
\Sh_n(y;p) =
\begin{cases}
\Sh_n(y;k_1-k_2+\frac{1}{2},k_1+k_2-\frac{1}{2},k_2-k_1-p+\frac{1}{2}),&
p \leq 0,\\
\Sh_n(y;k_2-k_1+\frac{1}{2},k_1+k_2-\frac{1}{2},k_1-k_2+p+\frac{1}{2}),&
p\geq0,
\end{cases}
\]
and
\[
d \mu(y;p) =
\begin{cases}
d
\mu(y;k_1-k_2+\frac{1}{2},k_1+k_2-\frac{1}{2},k_2-k_1-p+\frac{1}{2}),&
p \leq 0,\\ d \mu
(y;k_2-k_1+\frac{1}{2},k_1+k_2-\frac{1}{2},k_1-k_2+p+\frac{1}{2}),&
p\geq0.
\end{cases}
\]
\end{thm}
So the continuous dual Hahn polynomials have an interpretation as
Clebsch-Gordan coefficients for $\su(1,1)$. The inversion of
\eqref{decomp vec} can be given explicitly, e.g. for an element
\[
f \tensor e_{r-L} = \int_0^\infty f(x) e_{r-L} dx \in L^2(0,\infty)
\tensor \lt(\Z) \cong \dirint \lt(\Z) dx
\]
in the representation space of the direct integral, we have
\begin{equation} \label{inv decomp}
f \tensor e_{r-L} =
\begin{cases}
\displaystyle \sum_{p=0}^\infty (-1)^{p-r} \left[ \int_{\R} \Sh_p(y;r)
f(y) d\mu(y;r) \right] e_p \tensor e_{p-r}, & r \leq 0,\\ \displaystyle
\sum_{p=0}^\infty (-1)^{p} \left[\int_{\R} \Sh_p(y;r) f(y) d \mu(y;r)
\right] e_{p+r} \tensor e_p, & r \geq 0.
\end{cases}
\end{equation}
For the discrete components in theorem \ref{decomp} we can replace $f$
by a Dirac delta function at the appropriate points of the discrete
mass of $d \mu(\cdot;r)$.
\begin{proof}
First we work out the proof for $k_1-k_2 \geq -\frac{1}{2}$ and
$k_1+k_2 \geq \frac{1}{2}$, so $d\mu(y;n_1-n_2)$ only has a absolutely
continuous part. We denote the continuous part of the measure
$d\mu(y;a,b,c)$ by $W(y;a,b,c)$, so
\[
W(y;a,b,c) =  \left|
\frac{\Gamma(a+ix)\Gamma(b+ix)\Gamma(c+ix)}{\Gamma(2ix)} \right|
\frac{1}{\sqrt{2\pi\, \Gamma(a+b) \Gamma(a+c) \Gamma(b+c)}}, \quad
y=x^2.
\]
Observe that
\begin{equation} \label{W(c+1)}
\frac{W(y;a,b,c+1)}{W(y;a,b,c)} =\sqrt{ \frac{c^2+y}{(a+c)(b+c)}}.
\end{equation}

Define an operator $\La$ by gluing together all the operators $\La_p$,
$p \in \Z$,
\begin{equation} \label{Lambda}
\La( e_{n_1} \tensor e_{n_2} )=(-1)^{n_2} \int_{0}^\infty
\Sh_l(x^2;n_1-n_2) W(x^2;p) e_{n_1-n_2-L} dx ,
\end{equation}
where $l= \min\{n_1,n_2\}$ and $e_{n_1-n_2-L}$ inside the integral is
an orthonormal basisvector of the representation space for the
principal series representation $\pi^{x,\eps}$. Note that $\La$ is
unitary.

We use \eqref{pos} and \eqref{neg} to find the action of $B$ in the
tensor product for $n_1=n$ and $n_1 - n_2=p<0$;
\begin{align*}
\Lambda \Big(   &  \pi_{k_1}^+ \tensor \pi_{k_2}^- \big( \Delta(B)
\big) \, e_n \tensor e_{n-p} \Big)  \\ = & (-1)^{n-p}\int_0^\infty
\sqrt{(n+1)(2k_1+n)}\Sh_{n+1}(x^2;p+1) W(x^2;p+1) e_{p-L+1} d x \\ &-
(-1)^{n-p-1} \quad  \int_0^\infty
\sqrt{(n-p)(2k_2+n-p-1)}\Sh_n(x^2;p+1) W(x^2;p+1) e_{p-L+1} d x
\\ =&(-1)^{n-p} \int_0^\infty \sqrt{(p-L+k_1-k_2+L+\frac{1}{2})^2+x^2
}\  \Sh_n(x^2;p) W(x^2;p) e_{p-L+1} d x \\ =& \dirint \pi^{x,\eps}(B)
dx
 \circ \Lambda (  e_n \tensor e_{n-p} ),
\end{align*}
where $\pi^{x,\eps}$ denotes the principal series representation. Here
we use \eqref{W(c+1)} and the following contiguous relation
\begin{equation} \label{conti1}
\begin{split}
&\frac{c^2+y}{\sqrt{(a+c)(b+c)}}\, \Sh_n(y;a,b,c+1) = \\
&\sqrt{(n+1)(a+b+n)}\,\Sh_{n+1}(y;a,b,c) +
\sqrt{(a+c+n)(b+c+n)}\,\Sh_n(y;a,b,c).
\end{split}
\end{equation}
This relation can be verified by expanding
\[
(y+c^2)\Sh_{n}(y;a,b,c+1) = \sum_{i=0}^{n+1} c_i \Sh_i(y;a,b,c),
\]
where
\[
c_i= \int_\R (y+c^2) \Sh_n(y;a,b,c+1) \Sh_i(y;a,b,c) d\mu^2(y;a,b,c).
\]
From \eqref{W(c+1)} and the orthogonality relations for
$\Sh_n(y;a,b,c+1)$ we see that $c_i=0$ for $i<n$. The coefficients
$c_n$ and $c_{n+1}$ follow easily from the integral and the leading
coefficient of $\Sh_n(y;a,b,c) = ( \sqrt{ n!\, (a+b)_n (a+c)_n (b+c)_n
} )^{-1} y^n + \ldots$.

For the case $p \geq 0$ the intertwining property for the action of $B$
follows similarly using the contiguous relation
\begin{equation*}
\begin{split}
&\sqrt{(a+c-1)(b+c-1)}\Sh_n(y;a,b,c-1)=\\
&\sqrt{(n+a+c-1)(n+b+c-1)}\Sh_n(y;a,b,c) + \sqrt{n(n+a+b-1)}
\Sh_{n-1}(y;a,b,c) .
\end{split}
\end{equation*}
This companion formula to \eqref{conti1} is checked in the same way as
\eqref{conti1}.

So we now have shown
\begin{equation} \label{intertw}
\Lambda \circ \pi_{k_1}^+ \tensor \pi_{k_2}^- \big( \Delta(Y) \big)  =
{\sideset{}{^\oplus}\int\limits_{0}^\infty} \pi^{\rho,\eps}(Y) d\rho
\circ \Lambda ,
\end{equation}
for $Y=B$. The case  $Y=H$ of \eqref{intertw} is straightforward from
\eqref{pos}, \eqref{neg}, \eqref{princ} and \eqref{Lambda}. The case
$Y=C$ of \eqref{intertw} can be deduced as for $Y=B$, or it can be
deduced from the case $Y=B$ and the unitary properties of $\La$ and the
representations of $\su(1,1)$.

We now have proved the theorem for the case $k_1-k_2 \geq
-\frac{1}{2}$, $k_1+k_2 \geq \frac{1}{2}$. The other cases can be
proved in the same way, only some more bookkeeping is involved for the
discrete mass points.
\end{proof}

\begin{rem}
(i) For $k_1+k_2<\frac{1}{2}$ a complementary series representation
enters the decomposition. This fact has already been observed by \O
rsted and Zhang in \cite[\S 3]{OZ} for the case $k_1 = k_2$. They also
found continuous dual Hahn polynomials as Clebsch-Gordan coefficients.
For the group $SU(1,1)$ there is no complementary series representation
in the decomposition of the tensor product, since then $k_1, k_2 \in
\frac{1}{2}\N$. For the group case we refer to Kerimov and Verdiyev
\cite{KV}, but they do not mention that the Clebsch-Gordan coefficients
are continuous dual Hahn polynomials. For $k_1=k_2$ the tensor product
decomposition naturally occurs in the framework of canonical
representations, see Hille \cite{Hil} and references therein.

(ii) For $k_1-k_2 < -\frac{1}{2}$, the measure on the right hand side
of \eqref{decomp vec} has only a discrete part for $n_1-n_2 \leq K$.

(iii) The representation on the right hand side in theorem \ref{decomp}
are defined on the dense domain consisting of finite linear
combinations of elements as on the right hand side of \eqref{decomp
vec}. Then $\La$ is an intertwiner in the sense of \cite[\S 8.2]{Sch}.

\end{rem}

\section{Coupled and uncoupled eigenvectors of $X_c$}
In this section we consider a self-adjoint element $X_c$ in $\su(1,1)$.
We determine eigenvectors of this element in the four different
representations. Then we determine the Clebsch-Gordan coefficients for
this basis of eigenvectors in the tensor product. As a result we find a
general bilinear generating function for Meixner polynomials.

\subsection{Eigenvectors of $X_c$ in irreducible representations}
\label{sec 3} The Meixner polynomials (see \cite{AAR},\cite{KS}) are
defined by
\begin{equation} \label{Meixner}
M_n(x;\beta,c) = \F{2}{1}{-n,-x}{\beta}{1-\frac{1}{c}}.
\end{equation}
For $\beta>0$ and $0<c<1$, these are orthogonal polynomials with
respect to a positive measure on $\Z_{\geq 0}$. The orthonormal Meixner
polynomials
\[
\M_n(x) = \M_n(x;\beta,c) = \sqrt{\frac{(\beta)_n\, c^n}{n!}}
M_n(x;\beta,c)
\]
satisfy the following three-term recurrence relation
\begin{equation} \label{rec Mp}
\begin{split}
\frac{(c-1)(x+\frac{1}{2}\beta)}{\sqrt{c}} \M_n(x) = &
\sqrt{(n+1)(n+\beta)} \M_{n+1}(x) + \\ &-
\frac{(c+1)(n+\frac{1}{2}\beta)}{\sqrt{c}} \M_n(x) +
\sqrt{n(n-1+\beta)} \M_{n-1}.
\end{split}
\end{equation}
Define the weight function $w(x;\be,c)$ by
\[
w(x;\be,c) =  \frac{(\beta)_x}{x!}\, c^x (1-c)^{\beta},
\]
then the Meixner polynomials are orthonormal in $x$ with respect to
this weight function;
\[
\sum_{x=0}^{\infty} \M_n(x) \M_m(x) w(x;\be,c) = \delta_{mn}.
\]
From the definition it is clear that the Meixner polynomials are
self-dual, i.e.~invariant under interchanging $n$ and $x$. This gives
us another orthogonality relation
\begin{equation} \label{dual orth Mp}
\sum_{n=0}^{\infty} \M_n(x) \M_n(y) w(x;\be,c) = \delta_{xy}.
\end{equation}

The Meixner functions (see Masson and Repka \cite[\S 3]{MR} and Koelink
\cite[\S 4.4]{Koe}) are defined by
\begin{equation} \label{M function}
\begin{split}
\m_n(x)=\m_n (x;\lambda,\eps,c) =  \left( \frac{1-c}{\sqrt{c}}
\right)^{-n} & \frac{\sqrt{\Gamma(n+\lambda+ \eps+1) \Gamma(n+\eps -
\lambda)}}{\Gamma(n+1-x)}  \\ &\times \F{2}{1}{n+\eps+\lambda+1,
n+\eps-\lambda}{n+1-x}{\frac{c}{c-1}}.
\end{split}
\end{equation}
Here we use the analytic continuation to $\C \setminus [1,\infty)$ of
the $_2F_1$-series. The Meixner functions are orthonormal with respect
to a positive measure on $\Z$, if $0<c<1$, $0 \leq \eps <1 $ and
$-\frac{1}{2}
<\lambda < -\eps$ or $-\frac{1}{2} < \lambda < \eps-1$ or
$\lambda=-\frac{1}{2} + i \rho$, $\rho \geq 0$. If $n+\la +\eps+1$ or
$n+\eps-\la$ is a non-positive integer, the $_2F_1$-series terminates
and can be written as a Meixner polynomial. The Meixner funcions
satisfy the following three-term recurrence relation
\begin{equation} \label{rec Mf}
\begin{split}
\frac{(c-1)(\eps+x)}{\sqrt{c}}\m_{n}(x) = &\sqrt{(n+\eps+\lambda+1)
(n+\eps-\lambda)}\, \m_{n+1}(x) +
\\ &- \frac{(1+c)(n+\eps)}{\sqrt{c}}\, \m_n(x) + \sqrt{(n+\eps +
\lambda)( n+\eps-\lambda -1)}\, \m_{n-1}(x).
\end{split}
\end{equation}
Define the weight function $\w(x;\la,\eps,c)$ by
\[
\w(x;\la,\eps,c) = \frac{c^{-x}(1-c)^{-2\eps}}{\Gamma(\eps+x-\lambda)
\Gamma(\eps+\lambda+x+1)},
\]
then the Meixner functions are orthonormal with respect to this weight
function
\begin{equation} \label{ort Mf}
\sum_{x=-\infty}^\infty
 \m_m(x) \m_n(x) \w(x;\la,\eps,c) = \delta_{mn}.
\end{equation}
The Meixner functions also satisfy the dual orthonormality relation
\begin{equation} \label{dual Mf}
\sum_{n=-\infty}^\infty \m_n(x) \m_n(y) \w(x;\la,\eps,c) = \delta_{xy}.
\end{equation}
Note that for $\la=-\frac{1}{2}+i\rho$, the $_2F_1$-series in the
definition of the Meixner function is a Jacobi function
$\varphi_{2\rho}^{(n-x,n+x+2\eps)}\big(c/(1-c)\big)$, where the Jacobi
function is defined by
\begin{equation} \label{Jac fnct}
\varphi^{(\al,\be)}_\si(t) = \F{2}{1}{\frac{1}{2}(\alpha+\beta+1-i\si),
\frac{1}{2}(\alpha+\beta+1+i\si)} {\alpha+1}{-t}.
\end{equation}
For Jacobi functions, see \cite{Koo} or \cite[\S 7.4.3]{VK}.\\

We define the self-adjoint element ${X_c}$ in $\su(1,1)$ by
\[
X_c = -\frac{1+c}{2 \sqrt{c}} H + B -C.
\]
We first consider the action of $X_c$ in the positive and negative
discrete series representations.
\begin{prop} \label{ev +-}
For $x \in \Z_{\geq 0}$
\[
v^+_x=\sum_{n=0}^\infty \M_n(x;2k,c)\sqrt{w(x;2k,c)} \, e_n \in
\lt(\Z_{\geq 0}),
\]
is a unit eigenvector of $\pi_k^+(X_c)$ for the eigenvalue
$(c-1)(x+k)/\sqrt{c}$. And for $x \in \Z_{\geq 0}$
\[
v^-_x=\sum_{n=0}^\infty  \M_n(x;2k,c) \sqrt{w(x;2k,c)} \,
 e_n \in \lt(\Z_{\geq 0})
\]
is a unit eigenvector of $\pi_k^-(X_c)$ for the eigenvalue
$-(c-1)(x+k)/\sqrt{c}$. Moreover, $\{v_x^\bullet \}_{x\in \Z_{\geq
0}}$, $\bullet=+,-$, is an orthonormal basis of $\lt(\Z_{\geq 0})$.
\end{prop}
For the positive discrete series this is \cite[prop. 3.1]{KJ1}), and
for the negative discrete series this proceeds completely analogously.
Proposition \ref{ev +-} also shows that $\pi^+_k(X_c)$, respectively
$\pi_k^-(X_c)$, is essentially self-adjoint, see \cite{Akh}. The
self-adjoint extension is defined by the same formula on its maximal
domain, and the domain $\mathcal D$ of the self-adjoint extension can
be discribed in terms of the eigenvectors of proposition \ref{ev +-} as
\[
\mathcal D = \left\{ \sum_{x \in \Z_{\geq 0}} c_x v_x^+ \in
\lt(\Z_{\geq 0}) \Big| \sum_{x \in \Z_{\geq 0}} c_x (x+k) v_x^+ \in
\lt(\Z_{\geq 0}) \right\},
\]
and similarly for the self-adjoint extension of $\pi_k^-(X_c)$.

Next we consider the action of $X_c$ in the principal and complementary
series representations. We denote both these representations here by
$\pi^{\la, \eps}$, where $\la = -\frac{1}{2}+i\rho$, $\rho \geq 0$, in
the principal series.
\begin{prop} \label{ev pc}
For $x \in \Z$
\[
w^{\lambda,\eps}_x = \sum_{n=-\infty}^\infty \m_n(x;\lambda,\eps,c)
\sqrt{\w(x;\lambda,\eps,c)} \,e_n \in \lt(\Z)
\]
is a unit eigenvector of $\pi^{\lambda,\eps}(X_c)$ for the eigenvalue
$(c-1)(\eps+x)/\sqrt{c}$. Moreover, $\{w_x^{\la,\eps}\}_{x\in \Z}$ is
an orthonormal basis of $\lt(\Z)$.
\end{prop}
\begin{proof}
The proof follows using a doubly infinite Jacobi operator,
corresponding to \eqref{rec Mf} and its spectral decomposition given by
\eqref{ort Mf}, see \cite{Koe},\cite{MR}.
\end{proof}
\begin{rem} Note that the spectrum of $\pi^{\lambda,\eps}(X_c)$ is
independent of $\la$.
\end{rem}
We will denote the eigenvector $w^{\la,\eps}_x$ for the principal
series, by $v^{\rho,\eps}_x$, so $v_x^{\rho,\eps} =
w_x^{-\frac{1}{2}+i\rho,\eps}$.

Proposition \ref{ev pc} implies that $\pi^{\la,\eps}(X_c)$ is
essentially self-adjoint, see \cite{Koe}, \cite{MR}, and that its
self-adjoint extension is given by the same formula on its maximal
domain. In terms of the eigenvectors of proposition \ref{ev pc} the
domain is
\[
\mathcal D = \left\{ \sum_{x \in \Z} c_x w_x^{\la,\eps} \in \lt(\Z)
\Big| \sum_{x \in \Z} c_x (\eps+x) w_x^{\la,\eps} \in \lt(\Z) \right\}.
\]

\subsection{Eigenvectors of $X_c$ in the tensor product}
In this subsection we want to link the coupled and uncoupled
eigenvectors of $X_c$ in the tensor product $\pi^+_{k_1} \tensor
\pi^-_{k_2}$, and calculate the corresponding Clebsch-Gordan
coefficients. We discuss the method that reduces the problem to theorem
\ref{decomp} and then indicate another method, based on the Jacobi
transform, in remark \ref{rem Jac}.

We define elements $H_c$, $B_c$ and $C_c$ by
\begin{eqnarray*}
H_c & = & \phantom{-} \frac{1+c}{1-c}H-
\frac{2\sqrt{c}}{1-c}B + \frac{2\sqrt{c}}{1-c}C, \\
B_c & = & -\frac{\sqrt{c}}{1-c}H + \frac{1}{1-c} B -\frac{c}{1-c} C, \\
C_c & = & \phantom{-}\frac{\sqrt{c}}{1-c}H - \frac{c}{1-c} B +
\frac{1}{1-c} C.
\end{eqnarray*}
then $H_c=gHg^{-1}$, $B_c=gBg^{-1}$ and $C_c = gCg^{-1}$ with $g \in
SU(1,1)$, $g=\binom{a\ b}{b\ a}$, $a=1/\sqrt{1-c}$, $b=\sqrt{c/(1-c)}$.
Then the relations
\[
[H_c,B_c] = 2B_c, \quad [H_c,C_c] = -2C_c, \quad [B_c, C_c] = H_c,
\]
follow, and also
\[
H_c^* = H_c, \quad B_c^*=-C_c.
\]
Note that $H_c = \frac{2\sqrt{c}}{c-1} X_c$, so that $H_c$ act
diagonally  with respect to orthonormal bases in propositions \ref{ev
+-} and \ref{ev pc} in all representations under consideration.
Moreover, $\Om =-\frac{1}{4}( H_c^2 + 2H_c +4C_cB_c)$ by this argument,
or by a straightforward computation.

Next we want to calculate the actions of $\pi^+_k(B_c)$,
$\pi^-_k(B_c)$, $\pi^{\la,\eps}(B_c)$ on the eigenvectors of
$\pi^+_k(H_c)$, $\pi^-_k(H_c)$, $\pi^{\la,\eps}(H_c)$ of propositions
\ref{ev +-} and \ref{ev pc}. From $-4\Om =H_c^2 + 2H_c +4C_cB_c$,
$B_c^* = -C_c$, using the polar decomposition for the closure of
$\pi^+_k(B_c)$, we see that the domain of the closure of $\pi^+_k(B_c)$
equals the domain of the self-adjoint extension of $\pi^+_k(H_c)$. In
particular, $v_x^+$ is in the domain of the closure of $\pi^+_k(B_c)$,
and similarly for $\pi^+_k(C_c)$ and for the other representations.
From the commutation relations we get $\pi^+_k(B_c)v_x^+ = \al
v_{x+1}^+$. By proposition \ref{ev +-}
\[
\al \inprod{v_{x+1}^+}{e_0} = \inprod{\pi^+_k(B_c) v_x^+}{e_0} = -
\inprod{v_x^+}{\pi^+_k(C_c)e_0}
=
-\sqrt{w(x;2k,c)}\left(\frac{2k\sqrt{c}}{1-c}-
\frac{c\sqrt{2k}}{1-c}\M_1(x;2k,c) \right).
\]
From the recurrence relation for the Meixner polynomials, we find
$\M_1(x;2k,c) = \sqrt{2kc}-\frac{(1-c)x}{\sqrt{2kc}}$ and for the
weight function we have
\[
\frac{w(x;2k,c)}{w(x+1;2k,c)}=\frac{x+1}{(2k+x)c}.
\]
So we find
\[
\pi^+(B_c) v^+_x = -\sqrt{(x+1)(2k+x)}v^+_{x+1}.
\]
In the same way we can find the actions of $B_c$ and $C_c$ in the
discrete series. For the principal unitary and the complementary
series, we use the same method, but since there is no highest/lowest
weight vector in the representation space, we also need the following
contiguous relations
\begin{align*}
c(1-c)(1-z) \F{2}{1}{a,b}{c-1}{z} = &
(c-a)(c-b)z\F{2}{1}{a,b}{c+1}{z}\\ &+c\big(c-1-(2c-a-b-1)z \big)
\F{2}{1}{a,b}{c}{z},
\end{align*}
and
\[
\F{2}{1}{a,b}{c}{z} = \F{2}{1}{a,b}{c-1}{z}-\frac{abz}{c(c-1)}
\F{2}{1}{a+1,b+1}{c+1}{z}.
\]
For the first relation, see \cite[\S 2.8, eq.~30]{Erd}; the second
relation follows directly from
\[
(c-1)_n = \frac{c-1}{c-1+n}(c)_n.
\]
Then we have
\begin{align*}
\pi_k^+ (H_c) v^+_x & =  2(k+x) v^+_x,\\ \pi^+_k(B_c) v^+_x & =
-\sqrt{(x+1)(2k+x)}v^+_{x+1},\\ \pi^+_k(C_c) v^+_x & =
\sqrt{x(2k+x-1)}v^+_{x-1}  ,
\end{align*}
\begin{align*}
\pi_k^- (H_c) v^-_x & =  2(k+x) v^-_x,\\ \pi^-_k(B_c) v^-_x & =
\sqrt{x(2k+x-1)}v^-_{x-1},\\ \pi^-_k(C_c) v^-_x & = -\sqrt{(x+1)(2k+x)}
v^-_{x+1}  ,
\end{align*}
\begin{align*}
\pi^{\lambda,\eps} (H_c) w^{\lambda,\eps}_x & = 2(\eps+x)
w^{\lambda,\eps}_x, \\ \pi^{\lambda,\eps}(B_c) w^{\lambda,\eps}_x & =
\sqrt{(x+\eps-\lambda)(x+\eps+\lambda+1)}w^{\lambda,\eps}_{x+1},
\\ \pi^{\lambda,\eps}(C_c)w^{\lambda,\eps}_ x & = -
\sqrt{(x+\eps+\lambda)(x+\eps-\lambda-1)} w^{\lambda,\eps}_{x-1}.
\end{align*}
Comparing these actions to the actions of $H$, $B$ and $C$ on the
basisvectors $e_n$ for the representations in \eqref{pos}, \eqref{neg},
\eqref{princ} and \eqref{comp}, we see that the actions are equivalent,
where $x$ plays the role of $n$. Observe that, to get exactly the same
actions in the discrete series representations, we have to let $H_c$,
$B_c$ and $C_c$ act on $(-1)^x v_x$.
\begin{prop} \label{prop v+v-}
For $x_1,x_2 \in \Z_{\geq 0}$, $v^+_{x_1} \tensor v^-_{x_2}$ is an
eigenvector of $\pi^+_{k_1} \tensor \pi^-_{k_2} \big( \De(X_c) \big)$
for the eigenvalue $(c-1)(x_1-x_2+k_1-k_2)/\sqrt{c}$. And
\begin{equation} \label{v+v-}
\begin{split}
v^+_{x_1} \tensor v^-_{x_2} =& \sum_{n_1=0}^\infty
\sum_{n_2=0}^\infty \M_{n_1}(x_1;2k_1,c) \M_{n_2}(x_2;2k_2,c)
\sqrt{w(x_1;2k_1,c) w(x_2;2k_2,c)} \, e_{n_1} \tensor e_{n_2}
\\ =& (-1)^{x_1}\int_\R \Sh_x(y;x_1-x_2)
 v_{x_1-x_2-L}d\mu(y;x_1-x_2),
\end{split}
\end{equation}
where $x= \min\{x_1,x_2\}$ and $v_{x_1-x_2-L}$ is the eigenvector for
$X_c$ in the principal, complementary or negative discrete series
representation, depending on the value of $y$ in the same way as in
theorem \ref{decomp}.
\end{prop}
The second expression of \eqref{v+v-} states that the continuous dual
Hahn polynomials have an interpretation as Clebsch-Gordan coefficients
for the basis of eigenvectors of $X_c$.
\begin{proof}
That $v^+_{x_1} \tensor v^-_{x_2}$ is an eigenvector of $\pi^+_{k_1}
\tensor \pi^-_{k_2} \big( \De(H_c) \big)$ for the eigenvalue
$2(x_1-x_2+k_1-k_2)$ follows from proposition \ref{ev +-}. This also
gives the first expression for $v^+_{x_1} \tensor v^-_{x_2}$. We find
the second expression using \eqref{decomp vec} in theorem \ref{decomp}
and the observations above.
\end{proof}
\begin{rem} \label{rem Jac}
We indicate an alternative proof of proposition \ref{prop v+v-}
assuming, for simplicity, $k_1-k_2 \geq -\frac{1}{2}$ and $k_1+k_1 \leq
\frac{1}{2}$. From propositions \ref{ev +-}, \ref{ev pc} and theorem
\ref{decomp} we get that
\[
v_{x_1}^+ \tensor v_{x_2}^- = \int_0^\infty F(y;x_1,x_2)
v_{x_1-x_2-L}^{\sqrt{y},\eps} d\mu(y;x_1-x_2)
\]
for some function $F$ to be determined. Taking inner products with $e_0
\tensor e_0= \int_0^\infty e_{-L} d\mu(y;0)$ gives
\[
\begin{split}
\sqrt{ w(x_1;2k_1,c) w(x_2;2k_2,c) } = \int_0^\infty &F(y;x_1,x_2)
\m_{-L}(x_1-x_2-L;-\frac{1}{2}+i\sqrt{y},\eps,c)\\& \times
\sqrt{\w(x_1-x_2-L; -\frac{1}{2}+i\sqrt{y},\eps,c)} d\mu(y;x_1-x_2)
\end{split}
\]
and since $\m_{-L}(x_1-x_2-L;-\frac{1}{2}+i\sqrt{y},\eps,c)$ can be
expressed in terms of $\varphi_{2\sqrt{y}}^{(x_2-x_1,
x_1-x_2-2L+2\eps)}$, see \eqref{Jac fnct}, $F$ is the inverse Jacobi
transform \cite{Koo} of elementary functions. This Jacobi transform can
be calculated explicitly using a limit case of Koornwinder's formula
\cite[eq.~(9.4)]{Koo} stating that Jacobi polynomials are mapped onto
the Wilson polynomials. Using simple expansions and Saalsch\"utz's and
Kummer's summation formulas the explicit form for $F$ in terms of a
terminating $_3F_2$-series can be derived.
\end{rem}

\subsection{Bilinear generating functions for Meixner polynomials}
In this subsection we use the result of proposition \ref{prop v+v-} to
find a bilinear generating functions for Meixner polynomials.

We write the second expression of \eqref{v+v-} as
\begin{align*}
(-1)^{x_1}\int_\R \Sh_x(y;x_1-x_2) & v_{x_1-x_2-L}d \mu(y;x_1-x_2)=
\\ \sum_{r=-\infty}^\infty &\int_\R \Sh_x(y;x_1-x_2)
\m_r(x_1-x_2-L;-\frac{1}{2}+i\sqrt{y},\eps,c) \\ & \times
 \sqrt{\w(x_1-x_2-L;-\frac{1}{2}+i\sqrt{y},\eps,c)}\, e_r d\mu(y;x_1-x_2)
.
\end{align*}
To show that this also covers the case $k_1-k_2<-\frac{1}{2}$, where we
should have a finite sum of eigenvectors $v^-$, observe  that for
$y=-(k_1-k_2+\frac{1}{2}+j)=-(\eps-L+\frac{1}{2}+j)$, $0 \leq j \leq
K$, the $_2F_1$-series in the Meixner function terminates. Reversing
the sum
\[
\F{2}{1}{-n,b}{c}{x} = (-x)^n \frac{(b)_n}{(c)_n}
\F{2}{1}{-n,-c-n+1}{-b-n+1}{\frac{1}{x}}.
\]
then gives
\begin{align*}
\m_{n_1-n_2-L}&(x_1-x_2-L;-\frac{1}{2}+i\sqrt{y},k_1-k_2+L,c)
 \sqrt{\w(x_1-x_2-L;-\frac{1}{2}+i\sqrt{y},k_1-k_2+L,c)}  \\
 = &\
\M_{j-n_1+n_2}(j-x_1+x_2; 2k_2-2k_1-2j,c) \sqrt{w(j-x_1+x_2;
2k_2-2k_1-2j,c)}.
\end{align*}
So indeed we find a finite sum of eigenvectors $v_{j-x_1+x_2}^-$.

Now use \eqref{inv decomp} to find
\begin{align*}
\inprod{v^+_{x_1} \tensor v^-_{x_2} }{e_{n_1} \tensor e_{n_2}} =
\\ (-1)^{x_1+n_2} \int_\R &\Sh_{n}(y;n_1-n_2)
\Sh_x(y;x_1-x_2) \m_{n_1-n_2-L}
(x_1-x_2-L;-\frac{1}{2}+i\sqrt{y},\eps,c)
\\&\ \times \quad
\sqrt{\w(x_1-x_2-L;-\frac{1}{2}+i\sqrt{y},\eps,c)} d \nu(y;n_1-n_2,
x_1-x_2 ),
\end{align*}
where $d\nu(\cdot;n,x)$ is the restriction to the diagonal of the
product measure obtained from the measures $d\mu(\cdot;n)$ and
$d\mu(\cdot;x)$. From the first expression of \eqref{v+v-} we find
\[
\inprod{v^+_{x_1} \tensor v^-_{x_2} }{e_{n_1} \tensor e_{n_2}} =
\M_{n_1}(x_1;2k_1,c) \M_{n_2}(x_2;2k_2,c) \sqrt{w(x_1;2k_1,c)
w(x_2;2k_2,c)},
\]
so we have found the following formula
\begin{equation} \label{sum int}
\begin{split}
(-1)^{n_2+x_1} \int_\R &\Sh_{n}(y;n_1-n_2) \Sh_x(y;x_1-x_2)
\m_{n_1-n_2-L} (x_1-x_2-L;-\frac{1}{2}+i\sqrt{y},\eps,c)  \\ & \times\
\sqrt{\w(x_1-x_2-L;-\frac{1}{2}+i\sqrt{y},\eps,c)} d
\nu(y;n_1-n_2,x_1-x_2) \\ = &\
 \M_{n_1}(x_1;2k_1,c) \M_{n_2}(x_2;2k_2,c)
 \sqrt{w(x_1;2k_1,c) w(x_2;2k_2,c)}.
\end{split}
\end{equation}
This can be interpreted as a connection formula between two orthogonal
bases in $L^2(\Z_{\geq 0}^2)$.
\begin{thm} \label{sum thm}
The continuous dual Hahn and the Meixner polynomials, in the
notations \eqref{cont dHahn} and \eqref{Meixner}, satisfy, for
 $k_2 \geq k_1$,
\begin{align} \label{sum}
 S_{x_2}&(y;k_2-k_1+\frac{1}{2},k_1+k_2-\frac{1}{2},
k_1-k_2+x_1-x_2+\frac{1}{2}) \notag \\ \times& \
\frac{1}{\Gamma(p+x_2-x_1+1)}
 \F{2}{1}{p+k_1-k_2+\frac{1}{2}+i\sqrt{y},
p+k_1-k_2+\frac{1}{2}-i\sqrt{y}}{p+x_2-x_1+1}{\frac{c}{c-1}}
\\=&
\sum_{n=0}^\infty C_n S_n(y;k_2-k_1+\frac{1}{2},k_1+k_2-\frac{1}{2},
k_1-k_2+p+\frac{1}{2}) M_{n+p}(x_1;2k_1,c) M_{n}(x_2;2k_2,c), \nonumber
\end{align}
where
\[
C_n = (-1)^{x_1+x_2} \frac{ c^{n+x_1}(1-c)^{2k_1+p} (2k_1)_{x_1}
(2k_2)_{x_2} } {\Ga(n+1) \Ga(n+p+1)}.
\]
\end{thm}
\begin{proof}
First we concentrate on the case $x_1 \geq x_2$. Using the
orthonormality of $\Sh_n$ in \eqref{sum int} and writing the
polynomials in the usual normalization, we find
\begin{align}
 S_{x_2}&(y;k_2-k_1+\frac{1}{2},k_1+k_2-\frac{1}{2},
k_1-k_2+x_1-x_2+\frac{1}{2}) \notag \\ \times& \
\frac{1}{\Gamma(p+x_2-x_1+1)}
 \F{2}{1}{p+k_1-k_2+\frac{1}{2}+i\sqrt{y},
p+k_1-k_2+\frac{1}{2}-i\sqrt{y}}{p+x_2-x_1+1}{\frac{c}{c-1}} \notag
\\=&
\begin{cases}
\displaystyle \sum_{n=0}^\infty C_n^-
S_n(y;k_1-k_2+\frac{1}{2},k_1+k_2-\frac{1}{2}, k_2-k_1-p+\frac{1}{2})
&\\ \hspace{4cm} \times\  M_n(x_1;2k_1,c) M_{n-p}(x_2;2k_2,c), & p \leq
0,\\ \displaystyle \sum_{n=0}^\infty C_n^+
S_n(y;k_2-k_1+\frac{1}{2},k_1+k_2-\frac{1}{2}, k_1-k_2+p+\frac{1}{2})
&\\ \hspace{4cm} \times \ M_{n+p}(x_1;2k_1,c) M_{n}(x_2;2k_2,c), & p
\geq 0,
\end{cases} \label{expr1}
\end{align}
where
\begin{align*}
C_n^- = (-1)^{x_1+x_2+p} &\frac{c^{n-p+x_1} (1-c)^{2k_1+p} (2k_1)_{x_1}
(2k_2)_{x_2} } {\Ga(n+1) \Ga(n-p+1)} \\ &\times \quad \left| \frac{
\Ga(k_1-k_2+\frac{1}{2}+i\sqrt{y}) \Ga(k_2-k_1-p+\frac{1}{2}+i\sqrt{y})
} {\Ga(k_1-k_2+p+\frac{1}{2}+i\sqrt{y})
\Ga(k_2-k_1+\frac{1}{2}+i\sqrt{y})} \right|
\end{align*}
and
\[
C_n^+ = (-1)^{x_1+x_2} \frac{ c^{n+x_1}(1-c)^{2k_1+p} (2k_1)_{x_1}
(2k_2)_{x_2} } {\Ga(n+1) \Ga(n+p+1)}.
\]
Now consider the expression on the right hand side of \eqref{expr1} for
$p \leq 0$. This expression is also well defined for $p \geq 0$, since
then the factor $(\Ga(n-p+1))^{-1}$, which equals zero for $n < p$,
makes the summation start at $n=p$ instead of $n=0$. Using \eqref{cont
dHahn} and writing out the $_3F_2$-series in the continuous dual Hahn
polynomial as a sum, we get
\begin{align} \label{eq S}
S_n(y;&k_1-k_2+\frac{1}{2},k_1+k_2-\frac{1}{2}, k_2-k_1-p+\frac{1}{2})
\notag
\\ =& \frac{ \Ga(2k_1+n) \Ga(n-p+1) }{ \Ga(2k_1) } \sum_{j=p}^n
\frac{ (-n)_j (k_1-k_2+\frac{1}{2}+i\sqrt{y})_j
(k_1-k_2+\frac{1}{2}-i\sqrt{y})_j } { (2k_1)_j \Ga(1-p+j) (1)_j }
\notag \\ =& (-1)^p\, \left|
\frac{\Ga(k_1-k_2+p+\frac{1}{2}+i\sqrt{y})}
{\Ga(k_1-k_2+\frac{1}{2}+i\sqrt{y})}  \right|^2 \frac{ \Ga(2k_1+n)
\Ga(n+1)}{ \Ga(2k_1+p) \Ga(p+1)}  \\ & \times \sum_{l=0}^{n-p} \frac{
(p-n)_l (k_1-k_2+p+\frac{1}{2}+i\sqrt{y})_l
(k_1-k_2+p+\frac{1}{2}-i\sqrt{y})_l }{ (2k_1+p)_l (p+1)_l\, l! } \notag
\\ =& (-1)^p\, \left| \frac{\Ga(k_1-k_2+p+\frac{1}{2}+i\sqrt{y})}
{\Ga(k_1-k_2+\frac{1}{2}+i\sqrt{y})}  \right|^2
S_{n-p}(y;k_2-k_1+\frac{1}{2},k_1+k_2-\frac{1}{2},
k_1-k_2+p-\frac{1}{2}). \notag
\end{align}
Now we shift the index $n$ in the expression for $p\leq 0$ of
\eqref{expr1} to $n+p$, and we use the reflection formula for the
$\Ga$-function, to see that the expressions for $p \leq 0$ and $p \geq
0$ in \eqref{expr1} are the same.

For $x_1 \leq x_2$ we get from \eqref{sum int} a formula similar to
\eqref{expr1}, with $S_{x_1}(y;k_1-k_2+\frac{1}{2},
k_1+k_2-\frac{1}{2}, k_2-k_1+x_2-x_1+\frac{1}{2})$ as continuous dual
Hahn polynomial on the left hand side and different $\Ga$-factors  in
the constants $C_n^+$ and $C_n^-$ on the right hand side. We use
\eqref{eq S} again to see that both expressions on the right hand side
are the same. Then we use \eqref{eq S} again with $p=x_1-x_2$ and the
reflection formula for the $\Ga$-function, to obtain exactly the same
formula as for $x_1 \geq x_2$.
\end{proof}
\begin{rem} Several special cases of theorem \ref{sum thm} are of special
interest, or reduce to well-known results.

(i) For $x_2=0$ the continuous dual Hahn polynomial on the left hand
side of \eqref{sum} and $M_n(x_2;2k_2,c)$ reduce to $1$. This gives the
expansion of a Jacobi function, as function of $y$, in terms of
continuous dual Hahn polynomials with Meixner polynomials as
coefficients.

(ii) For $x_1=x_2=0$ in \eqref{sum} we find a special case of the
generating function for continuous dual Hahn polynomials
\[
(1-t)^{-\gamma} \F{3}{2}{\gamma,a+ix,a-ix}{a+b,a+c}{\frac{t}{t-1}}
= \sum_{n=0}^\infty \frac{ (\gamma)_n S_n(x^2;a,b,c)} {(a+b)_n
(a+c)_n n!}t^n, \quad \gamma \text{ arbitrary.}
\]

(iii) For $p=0$, $k_1=k_2$ and values of $y$ such that the
$_3F_2$-series in the continuous dual Hahn polynomials becomes summable
by Saalsch\"utz theorem \cite[\S 2.1.5, eq.~30]{Erd}
\[
\F{3}{2}{-n,a,b}{c,1+a+b-c-n}{1} = \frac{(c-a)_n (c-b)_n } { (c)_n
(c-a-b)_n},
\]
we find the dual orthogonality relations \eqref{dual orth Mp} for the
Meixner polynomials, which are equivalent to the orthogonality
relations of the Meixner polynomials.

(iv) An expression somewhat similar to \eqref{sum}, but simpler in
structure, for the Al-Salam and Chihara polynomials can be found in
\cite[Thm.~4.3]{IS}
\end{rem}

\begin{rem}
The element $X_c$, which we considered in section \ref{sec 3}, is a
special case of the element $X_a = -2a H +B-C \in \su(1,1)$, namely the
case for $|a|>1$. For $|a|<1$ and $|a|=1$ this element $X_a$ in the
discrete series representations has a continuous spectrum
$(-\infty,\infty)$ and $[0,\infty)$, which corresponds to
Meixner-Pollaczek and Laguerre polynomials respectively (see Koelink
and Van der Jeugt \cite[\S 3.1]{KJ1} and Van der Jeugt \cite[\S
2]{Jeu}). But in the principal and the complementary series, the
corresponding functions give much more difficulties than the Meixner
functions we used in this paper. For $|a|<1$ the spectral measure is
calculated in \cite[ex.~4.4.11]{Koe} and it turns out that the spectral
projection is on a two-dimensional space of generalized eigenvectors.
So we find orthogonality relations involving two linearly independent
"Meixner-Pollaczek"-functions, instead of just one function as for
$|a|>1$. This makes it much harder to compute the Clebsch-Gordan
coefficients for the bases of (generalized) eigenvectors of $X_a$ in
this way. However Hjalmar Rosengren pointed us that the calculations of
the Clebsch-Gordan coefficients might be possible from proving the
analogue of \eqref{sum} by brute force. In particular, Rosengren can
prove \eqref{sum} by direct methods.
\end{rem}

\section{Holomorphic and anti-holomorphic functions}
In this section we choose explicit realizations of $\su(1,1)$ and the
basisvectors $e_n$. Using these in \eqref{v+v-} we find a bilinear
generating function for the Meixner functions.\\

Consider the Hilbert space of holomorphic functions on the unit disk
$D$, with finite norm with respect to the inner product
\[
\inprod{f}{g} = \sum_{n=0}^\infty \frac{n!}{(2k)_n} f_n \overline{g_n},
\quad f(z) = \sum_{n=0}^\infty f_n z^n.
\]
The standard orthonormal basisvectors in this space are
\begin{equation} \label{e+}
e_n=e_n^k(z) = \sqrt{ \frac{(2k)_n}{n!} } z^n.
\end{equation}
The realization of the positive discrete series representation is
\begin{align*}
\pi^+_k (H) &= 2z \frac{d}{dz} +2k, \\ \pi^+_k (B) &= z^2 \frac{d}{dz}
+2kz, \\ \pi^+_k (C) &= -\frac{d}{dz}.
\end{align*}
An eigenvector for $\pi^+_k(X_c)$ for the eigenvalue
$(c-1)(x+k)/\sqrt{c}$ is
\[
v^+_x = \sum_{n=0}^\infty \M_n(x;2k,c) e_n^+, \quad x \in
\mathbb{Z}_+.
\]
Using the orthonormal basisvectors above, this eigenvector becomes a
generating function for the orthonormal Meixner polynomials (see
\cite[p.349]{AAR} or \cite{KS})
\begin{equation} \label{gen func}
v^+_x(z) =\sum_{n=0}^\infty \sqrt{ \frac{(2k)_n}{n!}} \M_n(x;2k,c) z^n
= (1-\frac{z}{\sqrt{c}})^x (1-\sqrt{c}z)^{-x-2k},
\end{equation}
where $x \in \Z_{\geq 0}$. This is a holomorphic function for
$|z|<1/\sqrt{c}$. The right hand side can also be found by solving the
first order differential equation
\[
\pi_k^+(X_c)\, y(z) = \frac{(c-1)(x+k)}{\sqrt{c}} \, y(z), \quad
y(0)=1.
\]

Consider the Hilbert space of anti-holomorphic functions on the unit
disk $D$, with the same inner product as above, now using $f(z) =
\sum_{n=0}^\infty f_n \oz^n$. The standard orthonormal basisvectors in
this space are
\begin{equation} \label{e-}
e_n=e_n^k(\overline{z}) = \sqrt{ \frac{(2k)_n}{n!} } \overline{z}^n.
\end{equation}
The realization of the negative discrete series representation is
\begin{align*}
\pi^+_k (H) &= -2\oz \frac{d}{d\oz} -2k, \\ \pi^+_k (B) &=
-\frac{d}{d\oz}, \\ \pi^+_k (C) &= \oz^2 \frac{d}{d\oz} +2k\oz.
\end{align*}
An eigenvector for $\pi^-_k(X_c)$ for the eigenvalue
$(1-c)(x+k)/\sqrt{c}$ is
\begin{equation} \label{gen func2}
v^-_x = \sum_{n=0}^\infty \M_n(x;2k,c) e_n^-=
(1-\frac{\overline{z}}{\sqrt{c}})^x (1-\sqrt{c}\overline{z})^{-x-2k},
\end{equation}
where $x \in \mathbb{Z}_+$. This is an anti-holomorphic function for
$|z|<1/\sqrt{c}$.

Next we use \eqref{inv decomp} to realize $f \tensor e_{r-L}$ in the
tensor product representation as functions of $z$ and $\ow$. To make
calculations easier we will assume $k_1-k_2 \geq -\frac{1}{2}$ and
$k_1+k_2 \geq \frac{1}{2}$, so the orthogonality measure for the
continuous dual Hahn polynomials has only an absolutely continuous
part. The result remains valid in the general case.
\begin{prop} \label{f tensor e}
In the explicit realizations \eqref{e+} and \eqref{e-}, we have for
$|z\ow|<1$
\begin{align} \label{f tensor e1}
(f&\, \tensor\, e_{r-L}) (z,\ow)  = \nonumber \\ & \frac{z^r
(1-z\ow)^{-2k_2}} {\sqrt{2\pi\, \Ga(2k_1) \Ga(2k_2)}} \int_0^\infty
f(\rho^2)   \left| \frac{ \Ga(k_1-k_2+r+\frac{1}{2}+i\rho)
\Ga(k_1+k_2-\frac{1}{2}+i\rho) \Ga(k_2-k_1+\frac{1}{2}+i\rho) } {
\Ga(2i\rho) } \right| \nonumber \\ & \qquad \qquad \times
 \frac{1}{\Ga(1+r)} \F{2}{1}{k_2-k_1+\frac{1}{2}+i\rho,
k_2-k_1+\frac{1}{2}-i\rho}{1+r}{\frac{z\ow}{z\ow-1}} d\rho.
\end{align}
\end{prop}
\begin{proof}
Substitute \eqref{e+} and \eqref{e-} in \eqref{inv decomp}, interchange
summation and integration, use the generating function for continuous
dual Hahn polynomials \cite{KS}, \cite[p.~349]{AAR}
\[
\sum_{n=0}^\infty \frac{S_n(x^2;a,b,c)}{n!\,(a+c)_n} t^n =
(1-t)^{-b+ix} \F{2}{1}{a+ix,c+ix}{a+c}{t}, \qquad |t|<1,
\]
and the transformation \cite[\S 2.1.4, eq.~22]{Erd}
\[
\F{2}{1}{a,b}{c}{t} = (1-t)^{-a} \F{2}{1}{a,c-b}{c}{\frac{t}{t-1}},
\]
then we find \eqref{f tensor e1} for $r \geq 0$. Following the same
procedure gives for $r \leq 0$
\begin{align*}
(f\, \tensor\, e_{r-L}) &(z,\ow)  = \frac{(-1)^{r} \ow^{-r}
(1-z\ow)^{-2k_1} }{\sqrt{ 2 \pi \Ga(2k_1) \Ga(2k_2)}} \\  & \times
\int_0^\infty f(\rho^2) \left| \frac{\Ga(k_1-k_2+\frac{1}{2}+i\rho)
\Ga(k_1+k_2-\frac{1}{2}+i\rho) \Ga(k_2-k_1-r+\frac{1}{2}+i\rho)}
{\Ga(2i\rho)} \right| \\ & \qquad \times \frac{1}{\Ga(1-r)}
\F{2}{1}{k_1-k_2+\frac{1}{2}+i\rho, k_1-k_2+\frac{1}{2}-i\rho} {1-r}
{\frac{z\ow}{z\ow-1}} d\rho.
\end{align*}
Use \cite[\S 2.1.4, eq.~23]{Erd}
\begin{equation} \label{tf1}
\F{2}{1}{a,b}{c}{t} = (1-t)^{c-a-b} \F{2}{1}{c-a,c-b}{c}{t},
\end{equation}
and for $r \in \Z$
\begin{equation} \label{Fr}
\begin{split} \frac{1}{\Ga(1-r)} \F{2}{1}{a,b}{1-r}{t} &= \sum_{k=r}^\infty
\frac{ (a)_k (b)_k} {\Ga(1-r+k) k!} t^k =  \\ \frac{(a)_r (b)_r
t^r}{r!} \sum_{p=0}^\infty \frac{(a+r)_p (b+r)_p }{(r+1)_p\, p!} t^p &
= \frac{\Ga(a+r) \Ga(b+r)} {\Ga(a) \Ga(b)} \frac{t^r}{\Ga(1+r)}
\F{2}{1}{a+r, b+r}{r+1}{t},
\end{split}
\end{equation}
and the reflection formula for the $\Ga$-function, to see that the
expressions for $r\leq 0$ and $r \geq0$ are the same.
\end{proof}
\begin{rem}
Note that the $_2F_1$-series inside the integral is a Jacobi function.
So we can consider $f \tensor e_n$ in this realization as the inverse
Jacobi transform of the function $f$ times $\Ga$-functions.
\end{rem}

Next we consider the eigenvectors in this realization. Using
proposition \ref{f tensor e} and interchanging summation and
integration, we find
\begin{align}
(f \tensor v_{x-L}) (z,\ow) =& \frac{(1-z\ow)^{-2k_2}} {\sqrt{2 \pi
\Ga(2k_1) \Ga(2k_2) }} \int_0^\infty f(\rho^2) \left|
\frac{\Ga(k_2-k_1+\frac{1}{2}+i\rho) \Ga(k_1+k_2-\frac{1}{2}+i\rho)}
{\Ga(2i\rho)} \right| \nonumber \\ \times& \sqrt{\w(x-L;\rho,\eps,c)}
\sum_{n=-\infty}^{\infty} \m_{n-L}(x-L;-\frac{1}{2}+i\rho,\eps,c)
\frac{ |\Ga(n+k_1-k_2+\frac{1}{2}+i\rho) |} {\Ga(1+n)}\nonumber \\
\times &
 \F{2}{1}{k_2-k_1+\frac{1}{2}+i\rho,
k_2-k_1+\frac{1}{2}-i\rho} {1+n} {\frac{z\ow}{z\ow-1}} z^n \, d\rho
\label{f tensor v1}
\end{align}
Note that the $_2F_1$-series has the same form, after applying
transformation \eqref{tf1}, as the $_2F_1$-series in definition of the
Meixner function. So we may consider the inner sum as a non-symmetric
Poisson kernel for the Meixner functions. To find conditions for the
convergence of the inner sum, we look at the asymptotic behaviour of
the Meixner function. The asymptotic behaviour has been determined in
\cite[\S4.4]{Koe} and \cite[\S3]{MR}. We have
\[
\m_n(x;\la,\eps,c) = C_1 n^{x+\eps-\frac{1}{2}} c^{n/2} \left(1+
\mathcal O \big( \frac{1}{n} \big) \right), \qquad n \rightarrow
\infty,
\]
where $C_1$ is a non-zero constant independent of $n$. The asymptotics
as $n \rightarrow -\infty$ follows from the connection formula
\cite[eq.3.8]{MR}, \cite[eq.4.4.7]{Koe} and observing that one term
vanishes for $x \in \Z$. This gives
\[
m_{-n}(x;\la,\eps,c) = C_2 n^{-x-\eps-\frac{1}{2}} c^{n/2} \left( 1 +
\mathcal O \big(\frac{1}{n} \big) \right), \qquad n \rightarrow \infty,
\]
where $C_2$ is a non-zero constant independent of $n$. Now we find that
the inner sum in \eqref{f tensor v1} converges for $|z|< 1/\sqrt{c}$
and $|w| < 1/\sqrt{c}$. So for these conditions we may indeed
interchange summation and integration for sufficiently smooth $f$, e.g.
$f \in \mathcal S$.

Using \eqref{gen func} and \eqref{gen func2}, we find an expression for
$v_p^+ \tensor v_{p-x}^-$, $x \leq 0$, in this realization
\begin{align*}
&(v^+_p \tensor v_{p-x}^-)(z,\ow) = v_p^+(z) \tensor v_{p-x}^-(\ow)
=\\&
 (1-z\sqrt{c})^{-2k_1}(1-\ow\sqrt{c})^{-2k_2} (1-c)^{k_1+k_2}
\sqrt{\frac{(2k_1)_p}{p!}} \left( \frac{\sqrt{c}-z}{1-z\sqrt{c}}
\right)^p \ \sqrt{ \frac{(2k_2)_{p-x}}{(p-x)!}} \left(
\frac{\sqrt{c}-\ow} {1-\ow\sqrt{c}} \right)^{p-x}.
\end{align*}
Observe that this has exactly the same structure as the realization for
$e_n \tensor e_{n-r}$. In the same way we can give an expression for
$v_{p+x}^+ \tensor v_{p}^-$, $x \geq 0$. Now we use the analogue of
\eqref{inv decomp} for the eigenvectors and the substitutions
\[
r \mapsto x, \qquad z \mapsto \frac{\sqrt{c}-z}{1-z\sqrt{c}}, \qquad
\ow \mapsto \frac{\sqrt{c}-\ow}{1-\ow\sqrt{c}},
\]
to find another expression for $f \tensor v_{x-L}$ in this realization,
\begin{equation} \label{f tensor v2}
\begin{split}
(f & \tensor v_{x-L})(z,\ow) = (-1)^x
(1-z\sqrt{c})^{-2k_1}(1-\ow\sqrt{c})^{-2k_2} (1-c)^{k_1+k_2}
\\ & \times \left( \frac{\sqrt{c}-z}{1-z\sqrt{c}} \right)^x \left( 1-
\frac{(\sqrt{c}-z)(\sqrt{c}-\ow)}{(1-z\sqrt{c})(1-\ow\sqrt{c})}
\right)^{-2k_2} \frac{1}{\sqrt{2\pi \Ga(2k_1) \Ga(2k_2)}} \\& \times
\int_0^\infty f(\rho^2) \left| \frac{ \Ga(k_2-k_1+\frac{1}{2}+i\rho)
\Ga(k_1+k_2-\frac{1}{2}+i\rho) \Ga(k_1-k_2+x+\frac{1}{2}+i\rho)}
{\Ga(2i\rho)} \right| \\ &\times \frac{1}{\Ga(1+x)}
\F{2}{1}{k_2-k_1+\frac{1}{2}+i\rho, k_2-k_1+\frac{1}{2}-i\rho} {1+x}
{\frac{ (\sqrt{c}-z)(\sqrt{c}-\ow)} {(1-c)(z\ow-1)}} d\rho.
\end{split}
\end{equation}
Now combining expression \eqref{f tensor v1} and \eqref{f tensor v2},
we find a non-symmetric type Poisson kernel for the Meixner functions.
Then substituting
\[
\sqrt{\frac{z}{\ow}} \mapsto t, \quad z\ow \mapsto s, \quad x-L \mapsto
x, \quad -L \mapsto y,
\]
(and still using $\eps=k_1-k_2+L$) we have
\begin{thm} \label{P kern}
For $x,y \in \Z$, $\sqrt{|cs|}<|t| < 1/\sqrt{|cs|}$, $0<|c|<1$ and
$0<|s|<1$ the Meixner functions, as defined in \eqref{M function},
satisfy
\begin{align*}
\sum_{n=-\infty}^{\infty} & \m_{n}(x;-\frac{1}{2}+i\rho,\eps,c)
\m_{n}(y;-\frac{1}{2}+i\rho,\eps,s)\, t^n = \\ &t^y
(1-s)^{2\eps+y}(1-c)^{2\eps+y} (1-\sqrt{cs}t)^{-2\eps-x-y}
(\sqrt{s}-\sqrt{c})^{x-y} c^{-\frac{x}{2}} s^{\frac{y}{2}}
\frac{|\Ga(x+\eps+\frac{1}{2}+i\rho)|^2} { \Ga(1+x-y)}
\\ & \times \F{2}{1}{-y-\eps+\frac{1}{2}+i\rho,
-y-\eps+\frac{1}{2}-i\rho} {1+x-y} {\frac{
(\sqrt{c}-\sqrt{s}t)(\sqrt{c}-\sqrt{s}/t)} {(1-c)(s-1)}}.
\end{align*}
\end{thm}
We can consider theorem \ref{P kern} as a non-symmetric type Poisson
kernel for the Meixner functions. For $\sqrt{c}$ and $\sqrt{s}$ we use
the principal value of the square root.
\begin{rem}
(i) Recall that the Meixner functions can be written in terms of Jacobi
functions, which contain the Bessel functions as a limit case
\cite[\S2.3]{Koo}. In this light we can consider theorem \ref{P kern}
as an extension of Graf's addition formula \cite[\S11.3]{Wat} for
Bessel functions of integer order.

(ii) The conditions on $t$ are exactly the same conditions for the
convergence of the generating functions for the Meixner polynomials
\eqref{gen func} and \eqref{gen func2}.

(iii) The $_2F_1$-series on the right hand side is the same
$_2F_1$-series as in the Meixner function
$\m_x(y;-\frac{1}{2}+i\rho,\eps,
\frac{(\sqrt{c}-\sqrt{s}t)(\sqrt{c}-\sqrt{s}/t)}
{(1-\sqrt{cs}/t)(1-\sqrt{cs}t)})$.

(iv) Theorem \ref{P kern} holds for all $\rho = \sqrt{y}$, $y \in \supp
\mu$, so the theorem also holds for the values of $\rho$ corresponding
to the complementary series representation. For the values of $\rho$
corresponding to the negative discrete series representations, $\rho =
i(\eps+\frac{1}{2}+j)$, the sum terminates at $n=j+1$. So this case
gives the non-symmetric Poisson kernel for the Meixner polynomials.

(v) If we put $t =1$ and $s = c \in \R$, we find the dual orthogonality
relations \eqref{dual Mf} for the Meixner functions.

(vi) From the substitution $-L \mapsto y$ it is not clear that $y \in
\Z$ instead of $-y \in \Z_{\geq 0}$. Using \eqref{Fr} and the
reflection formula for the $\Gamma$-function, we can show that the
right hand side of the summation formula is invariant under
interchanging both $c$ and $s$, and $x$ and $y$. And then since $x \in
\Z$, the same holds for $y$. This also shows that $c$ may be complex
valued.
\end{rem}

We may also try other realizations. Instead of the element $X_c$, we
may consider the element $X_c^\phi$, which is still self-adjoint,
defined by, cf. \cite{KJ2},
\[
X_c^\phi = - \frac{1+c}{2\sqrt{c}} H + e^{i\phi}B - e^{-i\phi} C, \quad
\phi \in [0,2\pi].
\]
And we realize the orthonormal basisvectors as Meixner polynomials, as
in \cite[\S 2]{KJ2}. This method leads to the summations formula
\eqref{sum} of theorem \ref{sum thm} for complex values of $c$.\\

\end{document}